\documentclass[11pt]{article}
\usepackage{amsmath, amssymb, graphics,epsfig}
\newtheorem{theorem}{Theorem}[section]

\newtheorem{proposition}[theorem]{Proposition}
\newtheorem{corollary}[theorem]{Corollary}

\newtheorem{problem}[theorem]{Problem}
\newtheorem{example}[theorem]{Example}
\newtheorem{remark}[theorem]{Remark}

\def\IC{{\mathbb C}}
\def\IR{{\mathbb R}}
\def\IN{{\mathbb N}}
\def\bA{{\bf A}}
\def\bB{{\bf B}}

\def\bu{{\bf u}}

\def\be{{\bf e}}
\def\b0{{\bf 0}}
\def\bx{{\bf x}}

\def\cB{{\cal B}}

\def\cE{{\cal E}}

\def\cL{{\cal L}}

\def\cK{{\cal K}}

\def\cS{{\cal S}}
\def\cT{{\cal T}}
\def\cH{{\cal H}}

\def\cB{{\cal B}}

\def\cS{{\cal S}}
\def\cT{{\cal T}}
\def\cH{{\cal H}}
\def\BH{{{\cal B}(\cH)}}
\def\BK{{{\cal B}(\cK)}}

\def\tr{{\rm tr}\,}

\def\span{{\rm span}\,}

\def\diag{{\rm diag}\,}
\def\conv{{\rm conv}\,}
\def\cl{{\bf cl}\,}
\def\Ra{{\Rightarrow}}
\def\[{\left [}
\def\]{\right ]}
\def\({\left (}
\def\){\right )}
\def\la{{\langle}}
\def\ra{{\rangle}}
\def\Ra{{\ \Rightarrow\ }}

\def\<{{\langle}}
\def\>{{\rangle}}
\def\1{{\bf 1}}

\def\bx{{\bf x}}
\def\by{{\bf y}}

\def\qed{\hfill\vbox{\hrule width 6 pt
\hbox{\vrule height 6 pt width 6 pt}}\medskip}
\newcommand{\bprob}{\begin{problem}}
\newcommand{\eprob}{\end{problem}}

\newcounter{example}
\begin{document}
%\openup 1\jot
\baselineskip 15.7pt

\title
{Numerical Range Inclusion, Dilation, and Operator Systems}

\author{Chi-Kwong Li and Yiu-Tung Poon}
\date{}
\maketitle

\begin{abstract}
Researchers have identified complex matrices $A$ 
such that a bounded linear operator $B$ acting on a Hilbert 
space will admit a dilation of the form $A \otimes I$ whenever  
the numerical  range inclusion relation  $W(B) \subseteq W(A)$ 
holds. Such an operator $A$ and the identity matrix will span a
maximal operator system, i.e., every unital positive map 
from $\span \{I, A, A^*\}$ to $\cB(\cH)$, the 
algebra of bounded linear operators
acting on a Hilbert space $\cH$, is completely positive.
In this paper, we identify $m$-tuple of matrices 
$\bA = (A_1, \dots, A_m)$
such that  any $m$-tuple of operators $\bB = (B_1, \dots, B_m)$ 
satisfying  the joint numerical range inclusion 
$W(\bB) \subseteq \conv W(\bA)$
will have a joint dilation of the form 
$(A_1\otimes I, \dots, A_m\otimes I)$.
Consequently, every unital positive map from 
$\span\{I, A_1, A_1^*, \dots, A_m, A_m^*\}$ to $\BH$ is completely positive.
New results and techniques are obtained relating to the study
of numerical range inclusion, dilation, and 
maximal operator systems.

\end{abstract}



\section{Introduction}

Let $\cB(\cH)$ be the set of bounded linear operators acting 
on the Hilbert space $\cH$ with inner product $\langle \bx,\by\rangle$. 
If $\cH$ has dimension $n$, we  identify $\cB(\cH)$ with
$M_n$ and  $\cH=\IC^n$ with the usual inner product $\langle \bx,\by\rangle = \by^*\bx$.  
The {\it numerical range of}  $A \in \cB(\cH)$ is defined and denoted by
$$W(A) = \{ \langle A\bx,\bx\rangle : \bx \in \cH, \langle \bx,\bx\rangle  = 1\}.$$
We say that an operator $B \in \cB(\cH)$ admits a dilation $A\in \cB(\cK)$
if there is a partial isometry $X: \cH \rightarrow \cK$ such that 
$X^*X = I_\cH$ and $X^*AX = B$. For simplicity, we will say that 
$B$ admits a dilation of the form $A \otimes I$ if    there is a Hilbert space $\cL$ such that $B$ admits a dilation of the form $A\otimes I_{\cL}  $.

Researchers have identified matrices $A$ such that
an operator $B \in \cB(\cH)$ has a dilation of the form  $A \otimes I$ whenever 
$W(B) \subseteq W(A)$; see \cite{An,Ar2,CL1,CL2,M,N}. We summarize the known result 
in the following.
 
\begin{theorem} \label{1.1} Let $A \in M_2$, or $A \in M_3$ unitarily similar to
$[a_0] \oplus A_1$ with $A_1 \in M_2$. Then a linear operator $B \in \cB(\cH)$ admits a dilation of the form $A\otimes I$ 
if and only if $W(B) \subseteq W(A)$.
\end{theorem}

It turns out that the above result can be reformulated in terms of maximal operator
systems. Recall that
an {\it operator system} $S$ of $\cB(\cH)$ is a self-adjoint subspace of  
$\cB(\cH)$ which  contains $I_{\cH}$. A linear map 
$\Phi: S\to \cB(\cK)$ is {\it unital} if  $\Phi\(I_{\cH}\)=I_{\cK}$,
$\Phi $ is {\it positive} if $\Phi(A)$ is positive semi-definite 
for every positive semi-definite  $A\in \cS$, and  
$\Phi$ is  {\it completely positive} 
if $I_k\otimes \Phi: M_k(\cS)\to M_k(\cB(\cK))$
define by  $(S_{ij}) \mapsto (\Phi(S_{ij}))$ 
is positive for every $k\ge 1$; e.g., see \cite{Choi} and \cite{Paulsen}
for some general background.

Suppose $A\in M_n$ and  $B\in \cB(\cH)$. Let $\cS = \span\{I_n,A,A^*\}$. 
Define a unital linear map $\Phi: \cS\to \cB(\cK)$  by 
$\Phi(aI + bA +cA^*) = aI + bB + cB^*$ for a given $B \in \BK$. 
By \cite[Lemma 4.1]{CL2}, $\Phi$ is positive 
if and only if $W(B)\subseteq W(A)$. 
On the other hand,  $\Phi$ is  completely positive if and only if $B$ 
has a dilation of the form $I\otimes A$; see Theorem \ref{2.2} in the next section.  
Therefore, Theorem \ref{1.1} can be restated as follows.

\begin{theorem} \label{1.2}
Suppose $A = A_0$ or $[a] \oplus A_0$ with $A_0 \in M_2$ and  
$B \in \cB(\cH)$. Then the map $\Phi: \span\{I, A, A^* \}\rightarrow 
\BK$ defined by
$$\Phi(aI + bA +cA^*) = aI + bB + cB^*, \qquad  a,b,c \in \IC.$$
is positive if and only if  $\Phi$ is completely positive.
\end{theorem}

Following the discussion in \cite[Theorem 3.22]{P},
an operator system $\cS$ is {\it a maximal
operator system} if  every unital positive map $\Phi: \cS \rightarrow
\BH$ is   completely positive. In such a case, 
we will say that $\cS$ is an OMAX. In particular, 
it was shown in \cite{P} that $\cS$ is an OMAX
if and only if for every positive integer $n$
a positive semi-definite operator 
operator $(b_{ij}) \in M_n(\cS)$ is the limit of 
a finite sum of operators of the form $S \otimes B$, where $S \in \cS$ and 
$B \in M_n$ are positive semi-definite operators. 
Despite this nice characterization, it is not easy  
to check or construct OMAX.

In this paper, we use the joint numerical range to study a maximal 
operator system $\cS \subseteq \BH$ of finite dimension.
The joint numerical range of $(A_1, \dots, A_m)\in \BH^m$ is defined by
$$W(A_1, \dots, A_m) = 
\{ (\la A_1x,x\ra, \dots, \la A_mx,x\ra):
x \in \cH, \la x,x\ra = 1\} \subseteq \IC^m,$$
here $\IC^m$ denotes the set of row vectors; we will also use $\IC^m$ 
to denote the set of column vectors.

We say that 
$(B_1, \dots, B_m) \in \cB(\cK)^m$ 
has a joint dilation $(A_1, \dots, A_m) \in \BH^m$ 
if there is a partial isometry $V: \cK \rightarrow \cH$
such that $V^*A_jV = B_j$ for all $j = 1, \dots, m$. If there is a Hilbert 
space $\cL$ such that $(B_1, \dots, B_m) $ has a joint dilation 
$(A_1\otimes I_{\cL}, \dots, A_m\otimes I_{\cL})\in \cB(\cH\otimes   \cL)$, 
we will simply say that  $(B_1, \dots, B_m) $ has a joint dilation 
$(A_1\otimes I , \dots,   A_m\otimes I)$.
With this definition, we have the following; see Theorem \ref{2.2}.

\begin{theorem}\label{1.3}
Let $\cS \subseteq M_n$ be an operator system with a basis
$\{I, A_1, \dots, A_m\}$ consisting of Hermitian matrices. Then
$\cS$ is a maximal operator system
if and only if every $(B_1, \dots, B_m) \in \cB(\cH)^m$ 
with $W(B_1, \dots,  B_k) \subseteq \conv W(A_1, \dots, A_m)$
has a joint dilation of the form $(A_1 \otimes I, \dots, A_m\otimes I  )$.
\end{theorem}

We will present some basic results in the next section
to facilitate the later discussion.
In Section 3, we will identify new maximal operator systems.
Our study provides new techniques and examples in the study of
numerical range inclusions, dilation and maximal operator systems.

It is easy to show that $(B_1, \dots, B_m)$ admits
a joint dilation of the form $(A_1 \otimes I, \dots, A_m\otimes I)$
is equivalent to $(B_1, \dots, B_m)$ admits a dilation 
of the form $(I \otimes A_1, \dots, I \otimes A_m)$.
We will use these two equivalent conditions in our discussion.

\section{Basic results}

We first summarize some basic results on the joint numerical range
$W(A_1, \dots, A_m)$ of $A_1, \dots, A_m \in \BH$; e.g., see \cite{LP} and its references. 
Since $A_j = H_j + iG_j$ with $(H_j,G_j) = (H_j^*,G_j^*)$ for 
$j = 1, \dots, m$, $W(A_1, \dots, A_m) \subseteq \IC^m$ can be 
identified with 
$W(H_1,G_1, \dots, H_m,G_m) \subseteq \IR^{2m}$. We can focus on the joint 
numerical range of self-adjoint operators.  Below are
some basic properties of the joint numerical range;
see \cite{LP} and its references.

\begin{theorem} \label{2.1} Let $T_1, \dots, T_m \in \cB(\cH)$ be self-adjoint operators.
\begin{enumerate}
\item[{\rm (a)}] The set $W(T_1, \dots, T_m)$ is bounded.
\item[{\rm (b)}] The set $W(T_1, \dots, T_m)$ is closed if 
$\dim \cH < \infty$.
Otherwise, it  may not be closed.
\item[{\rm (c)}] When $\dim \cH = 2$,  $W(T_1, \dots, T_m)$ 
is convex 
if and only if  $\dim \span\{I, T_1, \dots, T_m\} \le 3$.
\item[{\rm (d)}] Suppose $\dim \cH \ge 3$, and
$\dim \span\{I, T_1, \dots, T_m\} \le 3$.
Then $W(T_1, \dots, T_m)$ is convex.  
\item[{\rm (e)}] Suppose $\dim \cH \ge 3$ and 
$\dim \span\{I, T_1, \dots, T_m\} = 4$. 
Then there is a rank 2 orthgonal projection
$T_0$ such that $W(T_0, T_1, \dots, T_m)$ is not convex.
\end{enumerate}
\end{theorem}

Note that an  operator system 
$\cS \subseteq \BH$ always has a basis
$\{I, A_1, \dots, A_m\}$ consisting of self-adjoint
operators. The following is an extension of \cite[Lemma 4.1]{CL2}. 
The assertions are probably well known to researchers in the area, 
we include a proof here for completeness.

\begin{theorem}  \label{2.2} 
Let $\cS = \span \{I, A_1, \dots, A_m\} \subseteq \BH$
and $B_1, \dots, B_k \in \cB(\cK)$, where 
$A_1, \dots, A_m,$ $B_1, \dots, B_m$ are self-adjoint.
Define a linear map $\Phi: 
\cS \rightarrow \cT$  by 
$$\Phi(u_0 I +\mu_1 A_1 + \cdots + \mu A_m) 
= \mu_0 I +\mu_1 B_1 + \cdots +\mu_m B_m \quad 
\hbox{ for all } \mu_0 ,
\mu_1, \dots, \mu_m \in \IC.$$ 
\begin{itemize}
\item[{\rm (a)}] 
The map $\Phi$ is positive if and only if
$$W(B_1, \dots, B_m) \subseteq 
\cl ( \conv W(A_1, \dots, A_m)),$$
where $\cl(S)$ denotes the closure of  $S\subset \IR^n$.

\item[{\rm (b)}] 
If $(B_1 \dots, B_m)$ admits a dilation
of the form $(  A_1\otimes I , \dots,   A_m\otimes I)$, then $\Phi$ is 
completely positive.
If $\Phi$ is completely positive and $\dim \cH < \infty$, then
$(B_1, \dots, B_m)$ admits a dilation of the form 
$$(  A_1\otimes I , \dots,   A_m\otimes I).$$ 
\end{itemize}
\end{theorem}

\it Proof. \rm (a) Note that 
$(a_1, \dots, a_m) \in \conv(W(\bA))$ if and only if for any 
real vector $(u_1, \dots, u_m)$, 
$$u_0   +u_1a_1 + \cdots + u_m a_m \le \max \sigma(u_0 I +u_1A_1 + \cdots + u_m A_m).$$
Here $\sigma(H)$ denotes the spectrum of $H \in \BH$.
Thus, $W(B_1, \dots, B_m) \subseteq 
\conv(W(A_1, \dots, A_m))$
if and only if 
$u_0 I + u_1 B_1 + \cdots + u_m B_m \ge 0$
whenever the real vector $(u_0, \dots, u_m)$ satisfies 
$u_0 I + u_1 A_1 + \cdots + u_m A_m \ge 0$.
The assertion follows.

\medskip
(b) Suppose $(B_1,\dots,B_m)$ admits a dilation of the form $(  A_1\otimes I_{\cL} , \dots,   A_m\otimes I_{\cL})$
 for some Hilbert space $\cL$, then there exists 
$V: \cK \to   \cH \otimes \cL$,  such that $V^*V=I_{\cK}$ and 
$B_i= V^*A_iV$ for $i=1,\dots,m$. Therefore, $\Phi$ is completely positive.

Now, suppose $A_1, \dots, A_m \in M_n$  and  $\Phi: \cS\to \cB(\cK)$ is completely positive. 
Then by Arveson's Theorem \cite{Ar1}, 
$\Phi $ can be extended to $ \Phi: M_n \to \cB(\cK)$. 
By a result of Choi (see \cite{Choi} and 
\cite[Theorem 3.14]{Paulsen}), if $\{E_{11}, E_{12}, \dots, E_{nn}\}$
is the standard basis for $M_n$, then
$C = (\Phi(E_{ij}))\in M_n(\BK)$ 
is a positive operator. Let 
$C^{1/2} = [C_1 \dots C_n]$ so that 
$C_j: \cK \rightarrow \IC^n \otimes \cK$.
Because $\Phi$ is unital, if $V^* = [C_1^* \cdots C_n^*]$,
then  
$$I_\cK = \sum_{j=1}^n \Phi(E_{jj}) =  \sum_{j=1}^n C_j^* C_j
= V^*V.$$
Suppose $I_\cL = I_n \otimes I_\cK$.
Then for $\ell \in \{1, \dots, m\}$,
$$V^*( A_\ell \otimes  I_\cL)V = 
\sum_{i,j} (A_\ell)_{ij} (C_i^*C_j)
= \sum_{i,j} (A_\ell)_{ij} \Phi(E_{ij})
= \Phi(A_\ell)  = B_\ell.$$
Thus, $(B_1, \dots, B_m)$ admits a
joint  dilation of the form
$(A_1 \otimes I_\cL, \dots, A_m \otimes I_\cL)$.
\qed

\begin{remark} \label{rmk1} \rm
Note that the second statement in (b) 
may not hold if $\cH$ is infinite dimensional.
For example, if $A = \diag(1,1/2,\dots)$ and $B = \diag(0,1)$,
then $aI + bA \mapsto aI + bB$ is a unital completely positive map, but 
$B$ as no dilation of the form $A \otimes I$ because $0 \in W(B)$ and $0 \notin W(A\otimes I)$.
This example shows a subtle difference between the condition 
that $(B_1, \dots, B_m)$ has a dilation of the form 
$(A_1 \otimes I, \dots, A_m\otimes I)$ and the condition that
the unital positive map $\phi$ sending $A_j$ to $B_j$ for $j =1, \dots, m$
is completely positive. 
\iffalse
As we will see in Proposition \ref{4.4}, 
the operator system spanned by $\{I, A\}$ is indeed an OMAX.
\fi
\end{remark}

\medskip

Recall that $f: \IR^m \rightarrow \IR^m$ is an affine map if it
has the form $\bx \mapsto \bx R + \bx_0$ 
for a real matrix $R \in M_m$ and $\bx_0 \in \IR^m$, here $\IR^m$ denotes
the set of $1\times m$ real vectors.
The affine map is invertible if $R$ is invertible, and the 
inverse of $f$ has the form $y \mapsto \by R^{-1}- \bx_0 R^{-1}$.
One can extend the definition of affine map to an 
$m$-tuple of self-adjoint operators in $\BH$ by 
$$(A_1, \dots, A_m) \mapsto  (A_1, \dots, A_m) (r_{ij} I_{\cH})
+ (x_1  I_{\cH}, \dots, x_m  I_{\cH})$$
for a real matrix $R = (r_{ij}) \in M_m$ and  $(x_1, \dots, x_m)\in \IR^m$.
It turns out that real affine maps on $\IR^m$ and $\BH^m$ 
behave nicely in connection to positive maps, completely positive maps,
and the joint numerical range.
We have the following result which can be easily verified.

\begin{proposition}\label{prop3}
Let $\cS\subseteq \BH$ 
be an operator system with a basis
$\{I, A_1, \dots, A_m\}$, and $\Phi: \cS \rightarrow \BK$
a unital linear map defined by 
$\Phi(A_j) = B_j \in \BK$ for $j = 1, \dots, m$, where $A_1, \dots, A_m$, 
$B_1, \dots, B_m$
are self-adjoint.
Suppose $f$ is an invertible affine map such that 
$f(A_1, \dots, A_m) = (\tilde A_1, \dots, \tilde A_m)$ 
and $f(B_1, \dots, B_m) = (\tilde B_1, \dots, \tilde B_m)$.
\begin{itemize}
\item[{\rm (a)}]
Then $\Phi$ is positive (respectively, completely positive)
if and only if the unital map $\tilde \Phi$ defined by
$\tilde \Phi(\tilde A_j) = \tilde B_j$ for $j = 1, \dots, m$,
is positive (respectively, completely positive).
 \item[{\rm (b)}]
The $m$-tuple of operators $(B_1, \dots, B_m)$ 
has a joint dilation of the form
$(I \otimes A_1, \dots, I \otimes A_m)$ if and only if 
$(\tilde B_1, \dots, \tilde B_m)$ has a joint dilation of the form
$(I \otimes \tilde A_1, \dots, I \otimes \tilde A_m)$.
\item[{\rm (c)}]
For any real (unit) vector $(u_1, \dots, u_m)$
$$W(u_1B_1+ \cdots + u_m B_m) 
\subseteq  W(u_1A_1+ \cdots + u_m A_m)$$ 
if and only if for any real (unit) vector $(v_1, \dots, v_m)$
$$W(v_1 \tilde B_1 +  \cdots + v_m \tilde B_m) \subseteq 
 W(v_1\tilde A_1+ \cdots + v_m \tilde A_m).$$
\end{itemize}
\end{proposition}

Recall that a simplex in $\IR^m$ is a convex polyhedral set
with $m+1$ vertices.
 
\begin{theorem} \label{2.5}
Let $H_1, \dots, H_m \in \cB(\cH)$ be self-adjoint operators
such that 
$W(H_1, \dots, H_m)$ is a simplex in $\IR^{m}$. Then 
$(B_1, \dots, B_m) \subseteq \cB(\cK)^m$ has a joint dilation of the form 
$(I\otimes H_1, \dots, I\otimes H_m)$ whenever
$W(B_1, \dots, B_m) \subseteq W(H_1, \dots, H_m)$. In other words,
$\span\{I, H_1, \dots, H_m\}$ is a maximal operator system in $\cB(\cH)$.
\end{theorem}

\it Proof. \rm Suppose $\conv W(H_1, \dots, H_m)$ is a simplex. Then 
by the result in \cite{BL},
every vertex $(a_1, \dots, a_m)$ corresponding to a joint eigenvalue 
of $(H_1, \dots, H_m)$ such that $H_j x = a_j \bx$ for a unit vector $\bx$.
Thus, there is a unitary $U$ such that $U^*H_jU=[a_j] \oplus \tilde H_j$ for $j=1,\dots,m$. For simplicity, we will say that $(H_1, \dots, H_m)$ is untiarily similar to $\([a_1] \oplus \tilde H_1,\dots [a_m] \oplus \tilde H_m\)$.
Then 
$W(\tilde H_1, \dots, \tilde H_m)$ will contain  the other vertices of 
$W(H_1, \dots, H_m)$. We can then repeat the above argument, and extract
another joint eigenvalue $(b_1, \dots, b_m)$ of $H_1, \dots, H_m$.
Thus, $(H_1, \dots, H_m)$ is unitarily similar to $\(\diag(a_1, b_1) \oplus \hat H_1,\dots \diag(a_m, b_m) \oplus \hat H_m\)$ . Repeating this argument, we see that  $(H_1, \dots, H_m)$ is unitarily similar
$\( D_1 \oplus C_1,\dots , D_m \oplus C_m\)$ such that $D_j \in M_{m+1}$ is a diagonal matrix
and $W(D_1, \dots, D_m) = W(H_1,\dots, H_m) = \conv W(H_1, \dots, H_m)$.

Now, we show that $(B_1, \dots, B_m)$  admits a joint dilation of the form 
$(D_1 \otimes I, \dots, D_m\otimes I)$ for any self-adjoint operators
$B_1, \dots, B_m \in \cB(\cK)$ satisfying $W(B_1, \dots, B_m) \subseteq 
W(D_1, \dots, D_m)$. Our conclusion will follow.

By Theorem \ref{2.2}, we can apply an affine transform and assume that
$W(D_1,\dots, D_m)$ is the standard  simplex with vertices 
$\b0, \be_1, \dots, \be_m \in \IR^{m+1}$, where $\be_i$ has $1 $ at the $i^{\rm th}$ coordinate and $0 $ elsewhere.  Let
$D_j = E_{jj} \in M_{m+1}$ for $i=1,\dots,m$.
Then for every $k\ge 1$ and $C_0,\dots,C_m\in M_k$, 
$I \otimes C_0 + \sum_{j=1}^m D_j \otimes C_j \ge 0$
if and only if $C_0\ge 0$ and $C_0+C_j \ge 0$ for all $1\le j\le m$. By continuity argument, we may assume that $C_0$ is positive definite. Replacing $C_j$, with $C_0^{-1/2}C_jC_0^{-1/2}$, 
we may assume that $C_0=I$ and $C_j \ge -I$ for all $1\le j\le m$.   
If $W(B_1, \dots, B_m) \subseteq W(D_1, \dots, D_m)$, then
$B_j \ge 0$ for all $j$ and $\sum_{j=1}^m B_j \le I$.
Thus, 
$$I\otimes I + \sum_{j=1}^m B_j \otimes C_j
\ge I \otimes I + \sum_{j=1}^m B_j \otimes (-I) \ge 0.$$ 
\vskip -.5in \qed

\medskip
One may deduce \cite[Theorem 1.1]{BFL} from Theorem \ref{2.5} above.

\section{Maximal operator systems}

In this section, we identify some new maximal operator systems 
in addition  to the one described in Theorem \ref{2.5}.

\begin{theorem} \label{cs1cs2} Suppose $\cS_1 = \span\{I_A,A_1, \dots, A_r\}$
and $\cS_2 = \span\{I_B, B_1, \dots, B_s\}$ 
\iffalse
such that 
$$A_1, \dots, A_r, B_1, \dots, B_s$$
are invertible positive definite.
\fi
Then $\cS = \span (\{I_A\oplus 0, 0\oplus I_B\} \cup \{ A_i \oplus  0: 1\le i \le r\} \cup
\{0 \oplus  B_j: 1 \le j \le s\})$ is maximal 
if and only if 
$\cS_1$ and $\cS_2$ are maximal.
\end{theorem}

\it Proof. \rm Define $i_1:\cS_1\to \cS$, $i_2:\cS_2\to \cS$ , $\pi_1:\cS\to \cS_1$ and  
$\pi_2:\cS\to \cS_2$ by $i_1(A)=A\oplus 0$,
 $i_2(B)=0\oplus B$, $\pi_1(A\oplus B)=A$  and $\pi_2(A\oplus B)=B$. 
 
 Suppose $\cS_1$ and $\cS_2$ are maximal.  Given $\Phi:\cS\to \BH$ positive, let 
 $\Phi_j=\Phi\circ i_j$ for $j=1,2$. Then $\Phi_1$ and $\Phi_2$ are positive, hence, 
 completely positive. Therefore, 
 $\Phi=\Phi_1\circ\pi_1 +\Phi_2\circ\pi_2$ is also completely positive. This proves that $\cS$ 
 is maximal.
 
 Conversely, suppose $\cS$ is maximal. Given positive maps $\Phi_j:\cS_j\to \BH$, let 
 $\Phi=\Phi_1\circ\pi_1 +\Phi_2\circ\pi_2$. Then $\Phi $ is positive, hence, completely 
 positive. Therefore,  $\Phi_j=\Phi\circ i_j$, $j=1,2$ are also completely positive.
\qed
 
 By the above result, and the fact that 
 $\span\{E_{11}, E_{22}, E_{12}+E_{21}\} \subseteq M_2$ is an OMAX,
 see \cite{Choi,LP2019} and also Theorem \ref{main} below, we have the  following. 
 
\begin{theorem} \label{4.2}
Suppose $\cS$ is an operator system in $M_n$.
If up to a unitary similarity tranform, $\cS$ has a spanning set
which is a subset of $\{E_{jj}: 1 \le j \le n\} 
\cup \{E_{2j-1,2j}+E_{2j,2j-1}: 1 \le j\le n/2\}$,
then $\cS$ is an OMAX.
\end{theorem}

\iffalse

By Theorem \ref{4.2}, one easily deduce the following corollary, and then deduce 
Theorem \ref{2.5}.

\begin{corollary} If $\cS$ has a basis 
$\{I, A_1, \dots, A_m\}$  
consisting of Hermitian matrices such that
$W(\bA)$ is a simplex for $\bA = (A_1, \dots, A_m)$. 
Then $\cS$ is an OMAX.
\end{corollary}
\fi

\begin{corollary} If $\cS \subseteq M_3$ has a basis 
$\{I, A_1, A_2, A_3\}$ such that 
$W(A_1, A_2, A_3)$ is an ice-cream cone, i.e., the 
convex hull  of an ellipstical disk (a degenerated 
ellipsoid in $\IR^3$) and a point, then $\cS$ is an OMAX.
\end{corollary}

\it Proof. \rm Suppose $\cS$ satisfies the assumption. Then by an affine transform, we may 
assume that $W(A_1,A_2,A_3)$ is the ice-cream cone equal to the convex hull
of $\{(x,y,0): x^2+y^2 = 1\}$ and $\{(0,0,1)\}$ so that 
the matrices become $A_1 = E_{11} - E_{22}, A_2 = E_{12}-E_{21}, A_3 = E_{33}$.
Then $\span\{I, A_1, A_2, A_3\}$ has a basis $\{E_{11}, E_{22}, E_{33}, E_{12}+E_{21}\}$.
The result then follows from Theorem \ref{4.2}. 
\qed

Next, we focus on OMAX of the form $\span\{I, A, A^*\}$ for a
single matrix. Alternatively, we can write $A = A_1+iA_2$ for two Hermitian matrices $A_1$ and $A_2$ and consider
$\cS = \span\{I, A_1, A_2\}$. Clearly, if $\cS$ has dimension 
1, i.e., $A_1$ and $A_2$ are scalar operators,
then $\cS$ is an OMAX. We will consider the cases when
$\cS$ has dimension 2 and 3 in the following.

\begin{proposition} \label{4.4}
Suppose  $\cS \subseteq \BH$ has a basis $\{I, A\}$, with $A=A^*$.
Then $\cS$ is an OMAX. Furthermore, if $W(A)$ is closed, then
a bounded linear operator $B \in \BK$ has a dilation
of the form $A \otimes I$ whenever $W(B) \subseteq W(A)$.   
\end{proposition}

\it Proof. \rm We may replace $A$ by $\mu I + A$ and assume that 
$0 \in W(A)$ and $A$ has an operator matrix with the $(1,1)$ entry equal to zero.
\iffalse $A$ is positive semidefinite and $0\in \cl W(A)$. \fi
Suppose 
$\Phi:\cS \to \BK$ is a unital positive map and $B=\Phi(A)$. 
Then $W(B) \subseteq \cl  W(A)$ by Theorem \ref{2.2} (a). 

We will show that $\Phi$ is  completely positive.  
Suppose $k\ge 1$ and  $C_0, C_1 \in M_k$ is
such that  $I_\cH \otimes C_0 + A \otimes C_1$ is positive semidefinite. 
Since the $(1,1)$ entry of the operator matrix $A$ is assumed to be 0, we see that
the corresponding $(1,1)$ block of $I_\cH \otimes C_0 + A \otimes C_2$ equal to $C_0$ 
is positive semi-definite.
\iffalse
First, we are going to show that $C_0\ge 0$.  
Since $0\in \cl W(A)$, there is a sequence of unit vectors $\{\bx_n\}$ in $\cH$ such that $\la A\bx_n,\bx\ra \to 0$ as $n\to \infty$. For every unit vector $\by\in \IC^k$, we have 
$$(I_\cH \otimes C_0 + A \otimes C_1)(\bx_n\otimes \by),(\bx_n\otimes \by)\ra\ge 0 \ \mbox{ for all }n\ge 0.$$
As a result, 
$$\la C_0 \by,\by\ra +\la A\bx_n,\bx_n\ra \la C_1 \by,\by\ra\ge 0 \ \mbox{ for all }n\ge 0$$
so that 
$$\la C_0 \by,\by\ra \ge 0.$$
Therefore, $C_0\ge 0$. 
\fi
We may focus on the case when $C_0$ is positive definite, and obtain the conclusion 
by continuity argument. Replacing $C_j$ by $U^*C_0^{-1/2} C_j C_0^{1/2}U$
by a suitable unitary $U \in M_k$ for $j = 0,1$, we 
may assume that $C_0=I_k$ and $C_1=\diag(c_1,\dots, c_k)$ is a real diagonal matrix. 
Then for all  $1\le i\le k$,  we have 
$1+c_iA\ge 0$ implying $1+c_iB\ge 0$. Therefore, 
$I_{\cH}\otimes I_k+A\otimes  C_1\ge 0$. 
Since this is true for all $k \in \IN$, 
 $\Phi$ is completely positive.

The last assertion follows from Theorem \ref{2.5}.\qed

\baselineskip = 15pt

\begin{theorem}\label{main}
Suppose $\cS \subseteq M_n$ has dimension 3 and contains a rank 
one normal matrix. Then $\cS$ is an OMAX.
\end{theorem}

\it Proof. \rm %If $n\le 2$, then the result follows from Theorem \ref{1.2}. Suppose $n\ge 3$. 
By an affine transform we may assume that $A_1 = E_{11}$ and 
$A_2 = \begin{pmatrix} 0 & v^* \cr v& G\cr\end{pmatrix}$.
If $v = 0$, then $W(A) = \conv( \{1\} \cup W(iG))$
which is a line segment or a triangle. The result is known.

Suppose $v \ne 0$.
We can then replace $A_2$ by $([1]\oplus U^*)A_2([1]\oplus U)$ 
and assume that $v = (\gamma , 0, \dots, 0)^t$ with $\gamma > 0$.
We may replace $A_2$ by $A_2/\gamma$ and assume that $\gamma = 1$. Furthermore, by replacing 
$A_2$ with $A_2-aI_n-bA_1$,  we may assume that 
$G= \begin{pmatrix}0 &G_{12}  \cr G_{12}^*& G_{22}\cr\end{pmatrix}$, where $G_{22}\in M_{n-2}$.

Let $\Phi:\cS\to \BH$ be a unital positive map and $B_i=\Phi(A_i)$ for $i=1,2$. 
Since $A_1\ge 0$, we have $B_1\ge 0$. 
Suppose we have Hermitian matrices   $C_0,C_1, C_2 \in M_k$ such that 
$$I_n \otimes C_0 + A_1 \otimes C_1 + A_2 \otimes C_2=\begin{pmatrix}C_0+C_1&C_2&0 \cr
C_2&C_0&G_{12}\otimes C_2  \cr 
0&G_{12}^*\otimes C_2 & I_{n-2}\otimes C_0+G_{22}\otimes C_2 \cr\end{pmatrix}   \ge 0\,.$$
Therefore, $C_0$ is positive semidefinite. Without loss of generality, we may assume that 
$C_0=I_k$. We have
$$I_n \otimes I_k + A_1 \otimes C_1 + A_2 \otimes C_2 \ge 0$$
if and only if 
$$Q =   I_{n-1}\otimes I_k  + G \otimes C_2 \ge 0
\quad \hbox{ and } \quad I_k+C_1\ge 
(C_2 \ 0 \ \cdots \ 0)Q^\dag
(C_2 \ 0 \ \cdots \ 0)^*,
$$
where $X^\dag$ denotes the Moore-Penrose inverse of $X$.
For simplicity, 
we assume that the block matrix $Q$ is 
invertible and $C_2 = \diag(c_1, \dots, c_k)$.
Then we see that 
$$C_1 \ge  D = f(C_2),$$
where $D\in M_k$ is the diagonal matrix obtained by applying 
the  rational function 
$$f(x) = x^2 \det(I_{n-2}+xG_{22})\det(I_{n-1}+xG)^{-1} - 1\,.$$

Let  $D = \diag(d_1, \dots, d_k)$ and $C_1=D+P$ for some positive semidefinite $P\in M_k$. We have 

$$I_n  + d_iA_1   + c_iA_2   \ge 0\quad \mbox{ for all } 1\le i\le k\,.$$

Since $\Phi$ is positive, we have 

$$I_{\cH}  + d_iB_1   + c_iB_2   \ge 0\quad \mbox{ for all } 1\le i\le k\,.$$

Therefore, 
\begin{eqnarray*}
I_{\cH}\otimes I_k + B_1 \otimes C_1 + B_2 \otimes C_2
&=&I_{\cH}\otimes I_k + B_1 \otimes (D+P) + B_2 \otimes C_2\\
&\ge&I_{\cH} \otimes I_k+ B_1 \otimes D + B_2 \otimes C_2\\
&\ge&0\,.\end{eqnarray*}
\vskip -.3in \qed

\begin{corollary} 
If $A \in M_2$. Then $\span \{ I, A, A^*\}$ is an OMAX.

\noindent
If $A \in M_3$ and the boundary of $W(A)$ has a flat portion, 
then $\span\{I, A, A^*\}$ is an OMAX.
\end{corollary}

\it Proof. \rm If $A \in M_2$, then there is $A+A^* - aI$
is a rank one normal matrix for some $a \in \IR$. By Theorem 
\ref{main}, we get the conclusion.

If $A \in M_3$ and the bounary of $W(A)$ has a flat portion. 
Then we may replace $A$ by $e^{it}(A-\mu I)$ and assume that
$W(A) \subseteq \{x+iy: x \ge 0, y \in \IR\}$ and 
$W(A)$ contains a line segment joining $0$ to $ai$ for some 
$a > 0$. Then $A+A^*$ will be rank one.
By Theorem \ref{main}, we get the conclusion. \qed

Note that the above corollary covers Theorem \ref{1.2},
and identify some new matrices in $M_3$ such that
$\span\{I, A, A^*\}$ is an OMAX. For example, if 
$A = E_{11} + i G \in M_3$ for any Hermitian $G$,
then $\span \{I, A, A^*\}$ is an OMAX.

In fact,  if  $A  \in M_n$ with $n \ge 4$ and  $\span \{I, A, A^*\}$
contains a rank one normal  matrix,  then there is $a, b, c \in \IC$ such that
$a A +  b A^* + c I =  E_{11}$. Thus,  we may assume that $A = E_{11} + iG$ for
a Hermitian matrix $G$ with $(1,1)$ entry equal to 0. 
Let $\hat G$ be obtained from $G$ by deleting its 
first row and first column. If $\hat G$ is a scalar matrix $gI_{n-1}$, then 
$A$ is unitarily similar $A_0 \oplus gI_{n-2}$ with $A_0 = 
\begin{pmatrix} 1 & g_{12}i\cr \bar g_{12}i & g\,i\end{pmatrix}$ and $W(A) = W(A_0)$.
If $\hat G$ is not a scalar matrix, then the boundary of $W(A)$ contains 
a line segment $\conv\{(0,y):y\in\sigma(\hat G)\}$. 
However, even if the boundary of $W(A)$ has a line segment,
there does not seem to be 
an easy way to decide whether $\span \{I, A, A^*\}$ 
conatins a rank one normal matrix 
in terms of $W(A)$ if $A \in M_n$ with $n \ge 4$.
Nonetheless, when $n = 4$, we can use the above analysis to determine
whether $\span\{I, A, A^*\}$ in contains a rank one normal matrix
in terms of the structure of $W(A)$,
and identify another new family of 
$A$ such the $\span \{I, A, A^*\}$ is an OMAX.

\begin{corollary} \label{5.3}
Let $A \in M_4.$  Suppose $W(A)$ is the convex hull of an elliptical disk $\cE$
and two points  $ \alpha, \beta \in \IC \setminus \cE$   
such that the  line $L$ passing through $ \alpha$ and $ \beta$ is tangent to  $\cE$.
Then $\span\{I, A, A^*\}$ is an OMAX. 
\end{corollary} 

\it Proof. \rm Suppose $A\in M_4$ satisfies the hypothesis.
We may replace $A$ by $e^{it}(A-\mu I)$ and assume that   $L$ is the imaginary axis $\{ iy:   y \in \IR\}$ and
$0\in \cE \subseteq \{x+iy: x \ge 0, y \in \IR\}$. Let $\alpha=ai,\ \beta=bi$ for some $a,b\in \IR$.
Then $A$ is unitarily similar to $\diag(ai,bi) \oplus A_0$,
where $A_0 \in M_2$ with $W(A_0) \subseteq \{x+iy:
x\ge 0, y \in \IR\}$ and $0 \in W(A_0)$.
It follows that $A+A^*$ is a rank one matrix.
By Theorem \ref{main}, the result follows.
\qed

\begin{remark} \label{rmk2} \rm
Note that Corollary \ref{5.3}  also holds if we allow the elliptical disk $\cE$ to degenerate to a line segment   and $L$ intersects $\cE$ at an endpoint.  The proof also works with $\alpha =\beta$.  Therefore,  the corollary also covers Theorem \ref{1.2}. It also provides some new examples of OMAX. Let  
$A=\diag(1+i,1-i)\oplus \begin{pmatrix} 0 &2 \cr 0& 0\cr\end{pmatrix}$. Then $W(A) $ is the convex hull of the unit disk and $\{1+i,1-i\}$.  Therefore, $\span\{I, A, A^*\}$ is an OMAX. 
\end{remark}

\medskip

\iffalse

\begin{example} \label{112}
Let  $A = \diag(b_1i,b_2i) \oplus A_0$,   
where $b_1,b_2\in \IR $ and 
$A_0=[1]$ or $ \begin{pmatrix} 0 & i \cr i & 2 \cr
\end{pmatrix}$.   
Then $A+A^* $ has rank 1.  It follows from Theorem \ref{main} 
that  the operator system $\cS$ spanned by $\{I,A,A^*\}$ is OMAX.

We note that  $W(A)$ is the convex hull of 
$\{ b_1i, b_2i\}\cup W(A_0)$.  If $b_1>0>b_2$, 
The line segment joining $ b_1i$ and $ b_2i$ is tangent to 
the unit disk at $0$.   
The result in  Theorem \ref{1.2}   can be deduced from this 
example by applying an invertible affine transformation. 
\end{example}

If $A \in M_n$ with $n \ge 4$ and if $\span \{I, A, A^*\}$
contains a rank one matrix,  then 
the boundary of $W(A)$ always has a flat portion.
However, there does not seem to be 
an easy way to decide whether $\span \{I, A, A^*\}$ 
conatins a rank one normal matrix 
in terms of $W(A)$ if $A \in M_n$ with $n \ge 4$.

\fi

Next, we consider the case when $A$ is a direct sum of matrices in $M_1$ and $M_2$ and characterize those $A$'s for which 
$\span\{I, A, A^*\}$ is an OMAX.

 \iffalse
\medskip\noindent
{\bf Remark} By the above proposition,
there are examples of symmetric $A \in M_3$ 
such that $\span\{I, A,A^*\}$ is not an
OMAX.

\medskip
Consider the following.

\medskip\noindent
{\bf Example.} 
Let $R$ be the convex hull of the unit disk and $\{1+i,1-i\}$.
We can let $A = \diag(1+i,1-i) \oplus \begin{pmatrix} 0 & 0 \cr 2 & 0 \cr
\end{pmatrix}$. Then $W(B) \subseteq W(A)$ if and only if 
$B$ has a dilation of the form $I\otimes A$.

More generally, we have the following. 
\fi

\begin{theorem} \label{main2} Let $A \in M_n$ be a direct 
sum of matrices in $M_1$ and $M_2$. The following conditions
are equivalent.
\begin{itemize}
\item[{\rm (a)}] $\span\{I,A,A^*\}$ is a maximal operator system.
\item[{\rm (b)}]  
Every $B \in B(H)$ with   
$W(B) \subseteq W(A)$ admits a dilation of the form $A\otimes I $.
\item[{\rm (c)}]  
Every  $B \in M_2$ with  $W(B) \subseteq W(A)$  admits a dilation
of the form $A\otimes I$.
\item[{\rm (d)}]  
The matrix is unitarily similar to 
$A =\hat A_1 \oplus \hat A_2$ with 
$W(\hat A_2) \subseteq W(\hat A_1)$, and 
$\hat A_1$ satisfies one of the following:

\medskip
{\rm (d.1)} \ $\hat A_1 \in M_2$,

\medskip
{\rm (d.2)} \  $\hat A_1 = [a_1] \oplus A_0 \in M_3,$

\medskip
{\rm (d.3)} \ $\hat A_1 = \diag(a_1, a_2) \oplus A_0\in M_4$ such that
 $\conv\{a_1, a_2\}\cap W(A_0) = \{a_0\} \notin \{a_1,a_2\}$,
 
 i.e., the line segment joining $a_1$ and $a_2$ touches 
       a boundary point of $W(A_0)$ at $a_0 \notin \{ a_1, a_2\}$.\end{itemize}
       \end{theorem}

If we allow  $a_0 = a_2$ or $a_0 = a_1 = a_2$ in (d.3),
then (d.3) will cover the cases (d.2) and (d.1), respectively.

\it Proof. \rm The implication $(d) \Ra (a)$ follows from Theorem \ref{1.2} and Corollary
 \ref{5.3} that if $\hat A_1$ satisfies
(d.1)--(d.3), then an operator $B$ satisfies $W(B) \subseteq W(A) = W(\hat A_1)$ will have 
a dilation of the form $ \hat A\otimes I$, and hence a dilation of the form $A\otimes I$. 

The implications $(a)\Ra (b)\Ra (c)$ 
follows from definition. We are going to prove $(c)\Ra (d)$. 

Suppose    satisfies     $A=\oplus_{i=1}^mA_i\in M_n $  satisfies (c), where each $A_i$ is in $M_1$ or $M_2$ and 
irreducible.
Furthermore, we can assume that $W(A_i)\neq W(A_j)$ for 
$i\neq j$. 

We may assume that 
\begin{equation}\label{wkm}W(A_{k+1} \oplus \cdots \oplus A_m)
\subseteq W(A_1 \oplus \cdots \oplus A_k) = W(A)\neq W\(\oplus_{i=1}^{j-1}A_i \oplus_{i=j+1}^{k}A_i \) \
\end{equation}
 for all $1\le j\le k$.
By (\ref{wkm}),   the boundary of $W(A) $,  $\partial W(A) $, consists of  elliptic arcs and line segments. 
Consider the following cases for the boundary $\partial W(A):$

\medskip
\noindent
\begin{minipage}{0.73\textwidth}  
\baselineskip 14.5pt
{\bf Case 1} Suppose  $\partial W(A) $ contains two non-degenerate elliptic arcs $S_1$ and 
$S_2$ coming from  two summands, say, $A_1,A_2 \in M_2$. For $3\le i\le m$, $W(A_i)$ can only 
contain  a finite number of points  in $S_1\cup S_2$. 
Therefore, we can choose an exposed extreme 
point $\mu_i$ of  $W(A_i)$ for $i=1,2$ such that $\mu_i\not \in W(A_j)$ 
for $j>2$.  Consider the line segment joining $\mu_1$ and $\mu_2$. 
We can construct an elliptical disk $\cE$ with 
the line segment joining $\mu_1, \mu_2$ as major axis and minor axis of length
$d > 0$ with sufficiently small $d$ so that $\cE \subseteq W(A)$. 
Then there 
exists $B\in M_2$ such that $W(B) =  \cE$.
We are going to show that $B$ does not have a dilation to 
$A \otimes I$.
\end{minipage}

\vskip -1.8in

{\hfill\includegraphics[height=1.6in]{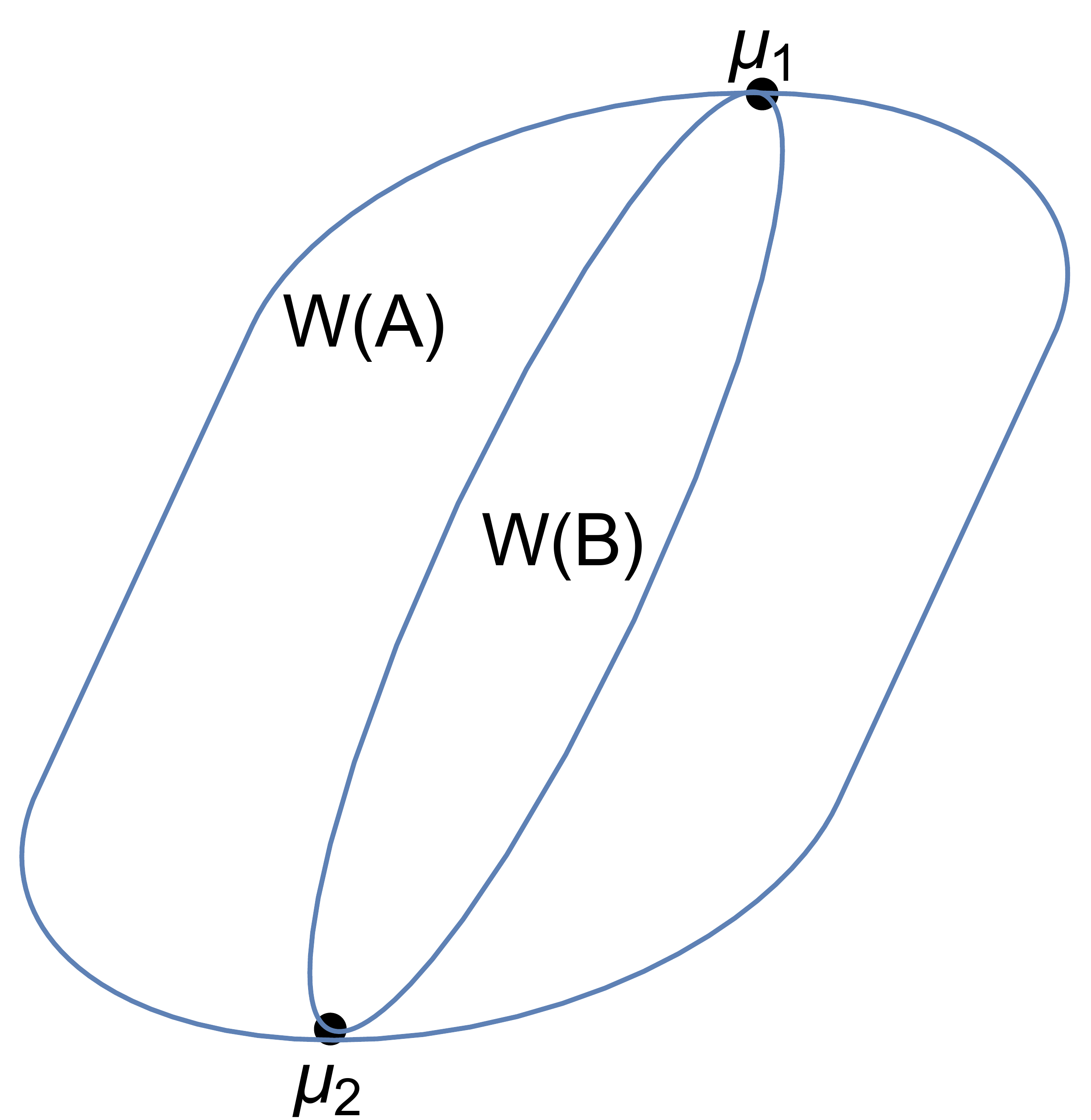}}

\newpage
\baselineskip = 16pt
Suppose the contrary that there exist $r\ge 1$ and $X\in M_{2\,rn}$ such that $XX^*=I_2$ and $X(A\otimes I_r)X^*=B$.  Let $\bu_1, \bu_2 $ be unit vectors such that $\mu_i=\bu_i^* B\bu_i $ for $i=1,2$. We may further assume that $\mu_i=(A_i)_{ii}$ for $i=1,2$. 

Let $X=[\bx_1\bx_2\cdots\bx_{rn}]$, where $\bx_j\in\IC^2$.  Since $\mu_1 $ (respectively, $\mu_2$) is an exposed extreme point of $W(B)$ and $W(A_1)$ (respectively, $W(A_2)$) , we have 

\begin{equation}\label{ux}\bu_1^*\bx_j=0\mbox{ for all }r<j\le nr\mbox{ and }\bu_2^*\bx_j=0\mbox{ for all }1\le j\le 3r\mbox{ and }4r<j\le nr\,.\end{equation} Since $\bu_1$ and $\bu_2$ are linearly independent, we have $\bx_j=0$ for all $r<j\le 3r$ and $4r< j\le nr$.
Also, 
$$\sum_{j=1}^r|\bu_1^*\bx_j|^2=\sum_{j=3r+1}^{4r}|\bu_2^*\bx_j|^2=1.$$
Thus, we have
\begin{eqnarray*}2
&=&\sum_{j=1}^r|\bu_1^*\bx_j|^2+\sum_{j=3r+1}^{4r}|\bu_2^*\bx_j|^2
\le\sum_{j=1}^r\|\bu_1\|^2\|\bx_j\|^2+
\sum_{j=3r+1}^{4r}\|\bu_2\|^2\|\bx_j\|^2\\
&\le&\sum_{j=1}^{nr}\bx_j^*\bx_j 
=\tr \(\sum_{j=1}^{nr}\bx_j\bx_j^*\) 
=\tr I_2 = 2.\end{eqnarray*}
Therefore, there exist $\alpha_j, \beta _j\in \IC$, $1\le j\le r$ such that  $\bx_j=\alpha_j\bu_1$ and and  $\bx_{3r+j}=\beta_j\bu_2$ for $1\le j\le r$. Hence, by (\ref{ux}), $\bu_1$ is orthoginal to $\bu_2$ and $B=\mu_1\bu_1\bu_1^*+\mu_2\bu_2\bu_2^*$ is normal, a contradiction.

\medskip
From the result in Case 1, $\partial W(A) $ can only contain  elliptic arcs from some $W(A_i)$ for  at most one $i$, with $1\le i\le k$. If $\partial W(A) $ does not contain any line segment, then condition   (d.1) is satisfied. It remains to consider the cases when 
$\partial W(A) $ contains some line segments.

\medskip
\noindent
\begin{minipage}{0.73\textwidth}
\baselineskip 15.5pt
{\bf Case 2}
Suppose  $\partial(W(A))$ has two pairs of consecutive line 
segments $\{L_1,L_2\} $ and $\{L_3,L_4\}$ with $L_1,L_2$ meeting at 
$\alpha$  and $ L_3,L_4 $ meeting at $\beta$ such 
that the open line segment $\overline{\alpha \beta}$ lies in the interior 
of $W(A)$. We may 
assume that $A=[\alpha]\oplus [\beta]\oplus_{i=3}^m A_i$ where 
$\alpha,\beta\not\in 
W(\oplus_{i=3}^m A_i)$. Let $A_0=\oplus_{i=3}^m A_i$.   
For $i=1,2$, let $p_i$ be the point on $L_i\cap W([\beta]\oplus A_0)$ nearest to $\alpha$.  
For  $i=3,4$, let $p_i$  be the point on $L_i\cap W([\alpha]\oplus A_0)$ nearest to $\beta$. 
We may apply an affine transform to $\IR^2$ and assume that
$(\alpha, \beta) = (-1,1)$, and 
$p_1, p_3, -p_1, -p_4$ have $y$-components larger than 2.  
We will show that  there is a circular disk $\cE$ in $W(A)$ with radius
less than 1 such that the boundary is tangent to  
at least 3 of the lines $L_i$'s, say $L_1,L_2,L_3$ at the points 
$\mu_1,\mu_2,\mu_3$ respectively. 
\end{minipage}

\vskip -2.1in 
{\hfill\includegraphics[height=2in]{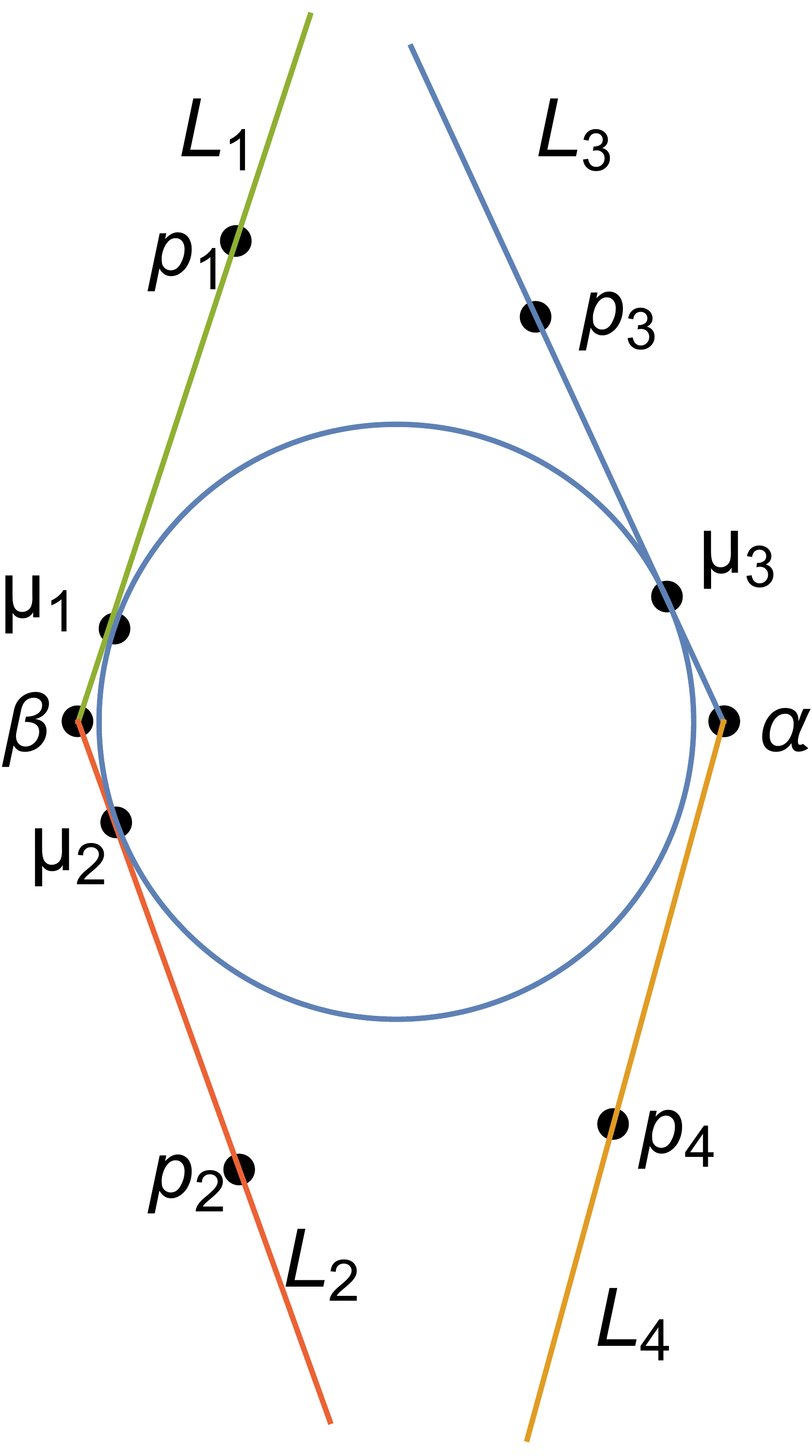} \ \ }

\bigskip
\rm 
Let $B_1$ and  $ B_2$
be the angular bisectors  at  $\beta$ and $\alpha$ respectively. 
\begin{itemize}
\item[{\rm 1)}]	If $B_1$ and $B_2$ coincide, then the midpoint $Q$ between $\alpha $ and $\beta$ has the same distance $r$  from $L_i$ for all $i=1,\dots,4$.
Let $\cE$ be the circle with center $Q$ and  radius $r$. 
 
\item[{\rm 2)}]	Suppose $B_1$ and $B_2$ intersect at a point $P$ and the distance from $P$ to $L_1$ (and $L_2$) is less than that from $P$ to $L_3$ (and $L_4$). If $B_1$ intersect $L_3$, then there exists $Q$ on $B_1$   such that the distance $r$ from $Q$ to $L_3$ is equal to that to  $ L_1$ (and $ L_2$), and less than the distance from $Q$   to $L_4$. Then $\cE$ is the circle with  $Q$ as the center and radius $r$. 
\end{itemize}

Now, suppose $B\in M_2$ with $W(B) = \cE$. 
We are going to show that $B$ does not have a dilation to $A\otimes I$ in the following.

Suppose $B=\sum_{j=1}^rX(j)AX(j)^{ *}$, 
where  $ \sum_{j=1}^rX(j) X(j)^{ *}=I_2$. Let 
$$X(j)=[\bx(j)_1\bx(j)_2\bx(j)_3] \ \hbox{ with } \
\bx(j)_1,\bx(j)_2\in \IC^2 \ \hbox{ and } \
\bx(j)_3\in M_{2\, {(n-2)}}.$$  If $\bu_1,\bu_2\in\IC^2$ are such 
that $\bu_i^*B\bu_i=\mu_i$ for $i=1,2$, then we must have 
$\bu_i^*\bx(j)_2=0$ for $i=1,2$ and all $1\le j\le r$. Since , 
$\bu_1$ and $\bu_2$ are linearly independent, we have $ \bx(j)_2=0$ for   all $1\le j\le r$. In this case, $\mu_3\not \in W(B)$, a contradiction. 

\medskip
From the results in Case 1 and 2, if $\partial W(A)$ only contains line 
segments, then it has  to be a (possibly degenerate) triangle, which is 
covered by (d.1) and (d.2). If $\partial W(A)$   contains an elliptic arc 
and some line segments. Then the line segments have to lie consecutively on 
$\partial W(A)$. Furthermore, there are either two or three
line segments. If $\partial W(A)$ contains an elliptic arc and  two line 
segments,  then $\hat A_1$ satisfies (d.2). So it remains to consider the 
case when $\partial W(A)$ contains an ellipic arc and three consecutive line 
segments.

\medskip

\noindent
{\bf Case 3} Suppose $\partial(W(A))$ contains three consecutive line segments $L_1,L_2,L_3$ 
and an elliptic arc $E$.  Suppose  $L_1$ and  $L_2$ meet  at $\alpha$, $L_2$ and  $L_3$ meet  
at $\beta$   and $E$ is part of the boundary of $W(A_i)$ for a summand $A_i\in M_2$ . 
If $L_2$ is tangent to the boundary of $W(A_3)$,  then, by   
Corollary \ref{5.3},  the operator system spanned by 
$\{I_4,[\alpha]\oplus[\beta]\oplus  A_3,[\overline{\alpha}]
\oplus[\overline{\beta}]\oplus  A_3^*\}$ is an OMAX. Consequently, 
the operator system spanned by $\{I_n, A,A^*\}$ is also OMAX. 

\medskip
\noindent
\begin{minipage}{0.6\textwidth}

\baselineskip 16pt
Suppose $L_2$ is not tangent to the boundary of $W(A_3)$. 
We may assume that $A_1=[\alpha]$,  $A_2=[\beta]$ and $i=3$. Therefore, $W(A)=W\([\alpha]\oplus[\beta]\oplus  A_3\)$. We may assume that $W(A_3)\not\subseteq W(A_i)$ for all $i>3$. Hence, there exists $\mu_3$ on $E$ such that $\mu_3\not \in W(A_i)$  for all $i\neq 3$.
We can construct $B\in M_2$ with $W(B)\subseteq W(A)$ and 
satisfies the following conditions: 
\baselineskip 13pt
\begin{enumerate}
\item  For  $i=1,2$, $L_i$ is tangent to $W(B)$ at $\mu_i$;
\item $W(B)$ and $W(A_3)$ have a common tangent at the point 
$\mu_3\in E$.
\end{enumerate}
\end{minipage}

\vskip -2.2in 
{\hfill\includegraphics[height=2in]{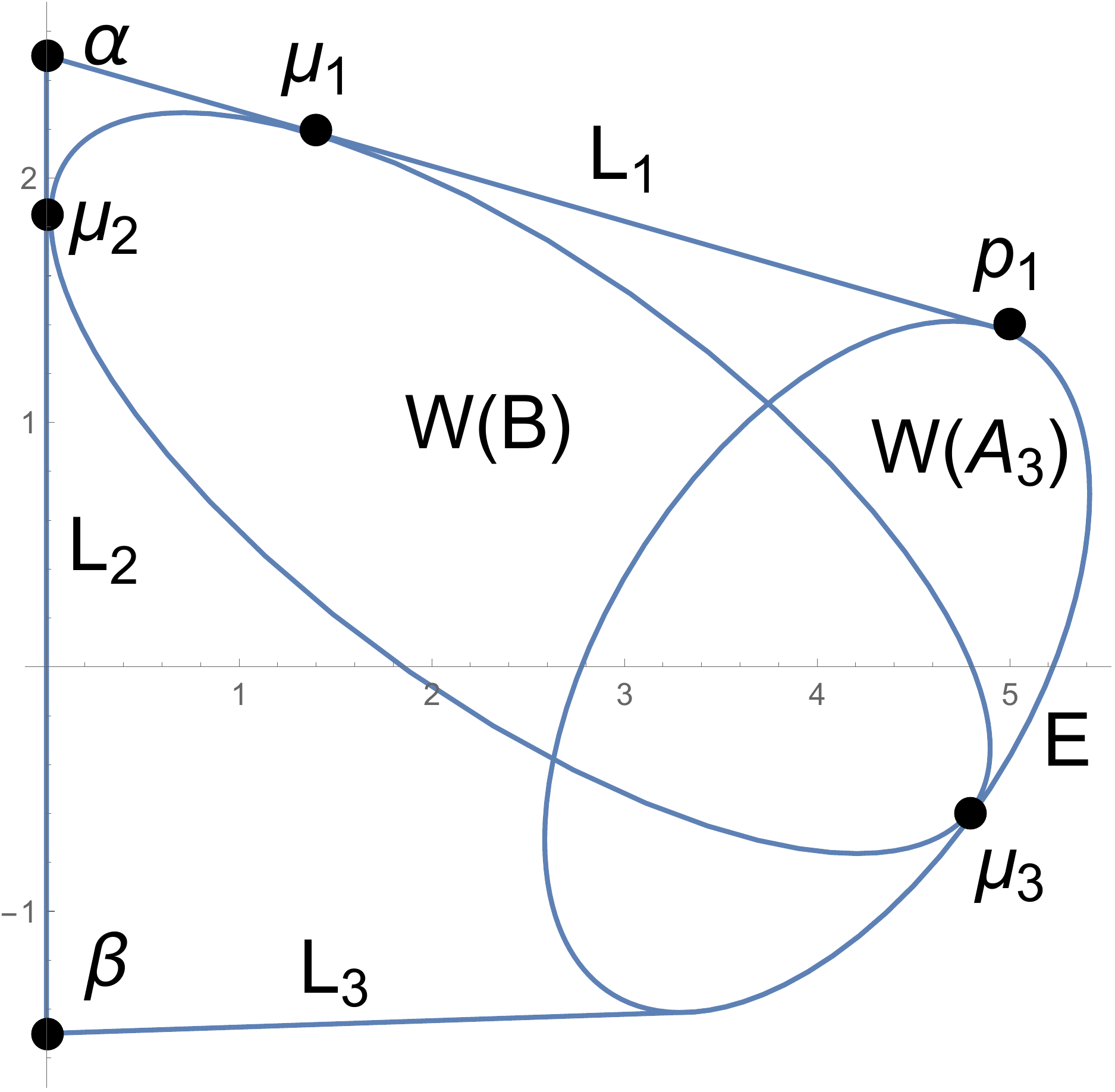} \ \ }

\vskip .2in\noindent
To see that such a matrix $B \in M_2$ exists, 
\iffalse
we can consider the triangle bounded  by the straight lines
$L_1, L_2$ and the tangent line at $\mu_3$.
Then there is
\fi
choose  a continuous family of ellipses $$\{\cE(\mu): \mu\mbox{ on the line segment joining $\alpha$ and $p_1$}\}$$ such that $L_1$ is a tangent to $\cE(\mu)$ at $\mu$  and $\mu_3$ is a boundary point of $ \cE(\mu)$ which  has a common tangent line with $W(A_3)$ at $\mu_3$. We may further assume that $\cE(p_1)\subseteq W(A_3)$.
\iffalse
there is a unique ellipse in the triangle with  $\mu_1, \mu_3$ as
as boundary points and touching the line $L_2$.
(Here we use the elementary fact that three tangent lines and two points 
determine a unique ellipse equation $x^2+axy + by^2 + cx + dy + f = 0$.)
If we choose $\mu_1$ sufficiently close to $\alpha$,
the ellipse  will lie in $W(A)$. Let $B \in M_2$ be such that 
the boundary of $W(B)$ is the ellipse. 
\fi
Since $\cE(p_1)\cap L_2=\emptyset $ and $\cE(\alpha)\cap L_2$ contains more than one points, there 
exists $\mu_1$ on the open  line segment joining $\alpha$ and $p_1$ such that $L_2$ is also 
tangent to $\cE(\mu)$.  Let $B \in M_2$ with boundary equal to $\cE(\mu_1)$.
Then $B$ will satisfy conditions (1) and (2) above.

We are going to prove that $B$ does not have a dilation to $A\otimes I$.
 Apply an affine map, if necessary, we can assume that 
$\alpha=b_1i$, $\beta =b_2i$ and $\mu_2= ci$ with $b_1>c>b_2$, and $W(A_3)$ lies on the open right half plane. 
Applying a unitary similarity to $B$, we can assume that 
$B=\diag(0, b)+iG$ for some $b>0$ and Hermitian $G$.  Then $\mu_2=e_1^*Be_1$,  where $e_1=\begin{pmatrix} 1\cr 0\end{pmatrix}$.

Suppose the contrary that $B=\sum_{j=1}^rX_jAX_j^*$ for some $X_j\in M_{n\,2}$ satisfying $\sum_{j=1}^rX_j X_j^*=I_2$. 
Let $X_j= [ x_1(j)  x_2(j) \cdots 
x_m(j)]$, where $x_i(j)\in \IC_{2}$.  Applying a unitary similarity to $A_3$, we can further assume that $A_3=\begin{pmatrix} *&*\cr *&\mu_3\end{pmatrix}$. 
For $i=1,3$, let $u_i\in \IC^2$ be a unit vector  satisfying $u_i^*Bu_i=\mu_i$.   For all $1\le j\le r$, we have 

\iffalse
\begin{enumerate}
\item  $u_1*x_2(j)=u_1*x_4(j)=0$;
\item $e_1^* x_3(j)=e_1^* x_4(j)=0$;
\item  $u_3^*x_1(j)=u_3^*x_2(j)=u_3^*x_4(j)=0$.
  \end{enumerate}
  \fi
  
 $$ u_1^*x_4(j)=0,\ 
 e_1^* x_3(j)=e_1^* x_4(j)=0 \mbox{ and }
 =u_3^*x_4(j)=0.$$

  Since $\{u_1,e_1\}$ and $\{e_1,u_3\}$ are linearly independent sets, we have $x_3(j)=x_4(j)=0$ for all $1\le j\le r$. Then $\mu_3\not \in W(B)$, a contradiction.
\qed

We restate the above result in terms of the geometrical shape of $W(A)$ 
in the following.

\begin{theorem} Suppose $A\in M_n$ is a direct sum of  matrices in $M_1$ and $M_2$. Then the operator system spanned by $\{I_n,A,A^*\}$
is an OMAX if and only if $W(A)$ is 
a singleton, a line segment, a triangular disk,
an elliptical disk, the convex hull of an elliptical disk $\cE$ 
with a point $\mu \notin\cE$, or 
the convex hull of $\cL$ and $\cE$, where 
$\cE$ is an elliptical disk $\cL = [z_1, z_2]$ 
is a line segment touching the 
ellipse $\cE$.  
\end{theorem}

\medskip\noindent
{\large \bf Acknowledgment}

Li is an affiliate member of the Institute for Quantum
Computing, University of Waterloo; he is also an honorary professor of Shanghai University. His
research was supported by USA NSF grant DMS 1331021, Simons Foundation Grant 351047, and
NNSF of China Grant 11971294.

\noindent 
(Li) Department of Mathematics, College of William \& Mary,
Williamsburg, VA 23187, USA.

\quad{ckli@math.wm.edu}

\noindent
(Poon)
Department of Mathematics, Iowa State University, Ames, IA 50011, USA.

\quad Center for Quantum Computing, Peng Cheng Laboratory, Shenzhen, 518055, China.
 
 \quad {ytpoon@iastate.edu}

\end{document}